# Non-homogeneous time convolutions, renewal processes and age dependent mean number of motorcar accidents

FULVIO GISMONDI^, JACQUES JANSSEN° and RAIMONDO MANCA*

^University "Guglielmo Marconi", Via Plinio, 44, 00193 Roma, Italy,

°Honorary professor at the Solvay Business School ULB, 4 Rue de la Colline, Tubize, Belgium 1480,

*University of Rome "La Sapienza", via del Castro Laurenziano 9, 00161 Roma, Italy,

raimondo.manca@uniroma1.it, tel. +390649766507, fax +390649766765.

'Work partially supported by a MIUR-PRIN and "La Sapienza" grants.

**Abstract** Non-homogeneous renewal processes are not yet well established. One of the tools necessary for studying these processes is the non-homogeneous time convolution.

Renewal theory has great relevance in general in economics and in particular in actuarial science, however most actuarial problems are connected with the age of the insured person. The introduction of non-homogeneity in the renewal processes brings actuarial applications closer to the real world. This paper will define the non-homogeneous time convolutions and try to give order to the non-homogeneous renewal processes. The numerical aspects of these processes are dealt with and, finally, a real data application to an aspect of motorcar insurance is proposed.

**Keywords** Non-homogeneous time, Convolution, Renewal processes, Numerical solution, Motorcar Insurance

**MSC 2000**: 60K05 (primary), 60-08, 91B30 (secondary)

## 1. Introduction

Non-homogeneous renewal theory has not been studied in depth; a check in Mathscinet reveals only two papers on the topic: (Wilson & Costello (2005) and Mode (1974)). On the other hand, a search for articles on non-homogeneous Poisson's processes shows 101 papers on the topic, but only 4 have relevance to actuarial science; two of them are written in Chinese (Cheng et al. (2006) and Chen & Yang (2007)) and the other two, (Lu & Garrido (2004) and Grandell (1971)) are in English. A search of actuarial literature in the IAA site reveals 95 papers and in the Casualty Actuarial Society site there are at least 10 papers dealing with non-homogeneous Poisson's processes. We recall Norberg (1993, 1999) who deals with the topic of claim reserve, as do Lu & Garrido (2004), Pfeifer & Nešlehovathá (2003) (who deals with the topic more in a financial environment), the introductory paper by Daniel (2008) and the book by Mikosch (2009). Considering the literature, it is clear that many researchers have worked with non-homogeneous Poisson's processes.

The non-homogeneous Poisson's process is a particular case of the non-homogeneous renewal process, but it has its "Achilles heel", i.e. the waiting time distribution function (d.f.), which rules that renewals must be negative exponentials. This constraint, as will be seen from our application, is usually not satisfied in real-world problems.

We think that studying the non-homogeneous renewal processes in a general way is necessary to work with non-homogeneous time convolutions. The papers Mode (1974) and Wilson & Costello (2004) take totally different approaches. The first paper established that the random variables (r.v.) of inter-arrival times were independent and not identically distributed, however the convolution operation among the non-homogeneous r.v. was not non-homogeneous. Regarding the second paper, the authors used an indirect approach, namely generalizing a method in a general renewal environment that was presented in Kallenberg (1975) for the definition of a non-homogeneous Poisson's process starting from a homogenous Poisson's process. This approach is not general because it supposes that the arrival times are independent Bernoulli r.v.. They used the thinning method to eliminate some of the arrival times. After the thinning application they assume that the inter-arrival times followed a gamma distribution. It is evident that, in a parametric environment, the probability distribution should be chosen by observed data and not decided in advance. Furthermore, a non-parametric study is not possible by means of this approach.

Regarding non-homogeneous convolutions again, it should be noted that references on this topic were found neither on Google nor in Mathscinet. It must also be pointed out that, despite the great relevance of this topic to actuarial science, no papers on general non-homogeneous renewal processes were found in the actuarial literature.

In a homogeneous renewal process, the system is renewed when the studied phenomenon is verified, then it restarts with the same initial characteristics. In the non-homogeneous case, the time when the renewal changes occurred changes the properties of the system. It is clear that a simple actuarial model, in which the time variable taken into account is the age, can be simulated well by this kind of stochastic process.

If the cumulating d.f. of the renewal process is of a negative exponential type, then the renewal process becomes a Poisson's process (see Janssen & Manca (2006) and Mikosch (2009)) and the related integral equation, also in the non-homogeneous case, can be solved analytically (Çynlar (1975)). This is a very particular case. In a more general case, the renewal equation cannot usually be solved analytically and it is necessary to solve the equation numerically.

In this light, this paper presents a straightforward method that deals directly with the numerical solution of the renewal equation. It proves that any kind of renewal process can be solved by means of this approach and that, in a very particular case of the general numerical solution (general because the d.f. is not specified), the related discrete time renewal equation can be obtained.

Freiberger & Grenander (1971) presented this approach in the homogeneous case without any theoretical justification, as did Xie (1989) in which the so-called "midpoint" formula was given with neither a general introduction nor a theoretical justification. Many other papers deal with the problem of discretization of the homogeneous renewal equation; see, for example, De Vylder & Marceau (1996) and Elkins & Wortman M.A. (2001), but as far as the authors are aware none of these papers state that it is possible to obtain the discrete time renewal equation from the continuous one and the continuous from the discrete. The first paper to state this relation in the homogeneous case in a Markov renewal environment is Corradi et al. (2004). In the book Janssen & Manca (2006), the same relation for the homogeneous renewal process case is presented for the first time. Regarding non-homogeneity, Janssen & Manca (2001) give these results in a non-homogeneous semi-Markov environment. More recently, Moura & Droguett (2008, 2010) presented a faster algorithm, which simplifies the problem and is useful for solving the non-homogeneous and the homogeneous semi-Markov evolution equation. Research aimed at numerically solving the non-homogeneous renewal equation has never been presented.

Once the strong relation between the continuous and discrete time renewal processes is proven, it is possible to work with a general discrete time d.f. In this light, it is possible to think of constructing the non-homogenous discrete time d.f. directly from the observed data. This approach leads to a d.f. that is constructed by the cumulative frequencies of the observed data. In this way, the model to be applied needs only raw data derived from observations.

The results derive directly from the data. The so-called "physical measure" is applied. The results are a direct function of the observed data.

All the theoretical apparatus will be applied to a motorcar insurance environment. By means of the non-homogeneous renewal equation, the mean number of accidents that an insured person can have within one year, two years, three years and so on, is computed. The non-homogeneity allows the production of these results for each considered age of insured people. We decided to consider the age as a non-homogeneous time variable because it is well known that drivers have different claims experiences in relation to their age. Indeed the insurance companies give different premiums dependent on age.

Applications of Poisson's processes and their generalizations concerning motorcar insurance can be found in Lemaire (1995). The author uses these processes in order to find the distribution of the number of accidents in one year.

The application of the renewal process given in our paper shows how to construct the mean number of accidents involving an insured person within any year of his driving life taking also into account the starting age. The d.f. put into the renewal function were constructed by real data provided by an insurance company.

The paper is organized as follows. In the second section the continuous time non-homogeneous convolutions and their properties are presented for the first time. The third section describes the continuous time non-homogeneous renewal processes and their evolution equation. In this section, also for purposes of completeness, the well-known properties of the homogeneous renewal process are reported. In the fourth section, the numerical solution of the evolution equation of the non-homogeneous continuous time renewal process is presented. Furthermore, how to obtain the discrete time evolution equation from the continuous time evolution equation and vice versa is explained. In section 5, an application to motorcar insurance in, both, the homogeneous and non-homogeneous environment is reported. Section 6 presents some short conclusive remarks.

**2. Continuous time non-homogeneous convolutions**

The following definitions are reported for reasons of clarity.

**Definition 2.1** A two variable function $f(s,t), 0 \leq s \leq t$, where $s$, $t$ represent times, is time non-homogeneous if:

$$\exists (s,t) \neq (s',t'), t - s = t' - s' : f(s,t) \neq f(s',t') \qquad \square$$

**Definition 2.2** Given two time non-homogeneous functions $f(s,t), g(s,t)$ their convolution is defined in the following way:

$$(f * g)(s,t) = \int_s^t g(s,\tau) f(\tau,t) d\tau; \quad s,t \in \hat{\mathbb{R}} = [-\infty, +\infty], s < t. \qquad \Box$$

Properties of non-homogeneous convolution:

**Associativity:** given $f(s,t), g(s,t), h(s,t)$) then:

$$(f * (g * h))(s,t) = ((f * g) * h)(s,t)$$

**Proof:**

$$(f * (g * h))(s,t) = f * \left( \int_s^t h(s,\tau_1) g(\tau_1,t) d\tau_1 \right) = \int_s^t \left( \int_s^{\tau_2} h(s,\tau_1) g(\tau_1,\tau_2) d\tau_1 \right) f(\tau_2,t) d\tau_2$$

$$= \int_s^t \int_s^{\tau_2} (h(s,\tau_1) g(\tau_1,\tau_2)) f(\tau_2,t) d\tau_1 d\tau_2 = \int_s^t \int_{\tau_1}^t h(s,\tau_1) g(\tau_1,\tau_2) f(\tau_2,t) d\tau_2 d\tau_1$$

$$= \int_s^t h(s,\tau_1) \int_{\tau_1}^t g(\tau_1,\tau_2) f(\tau_2,t) d\tau_2 d\tau_1 = \int_s^t h(s,\tau_1)(f * g(\tau_1,t)) d\tau_1 = ((f * g) * h)(s,t) \qquad \Box$$

**Remark 2.1.** The fourth step of the proof is the application of Dirichlet formula. $\qquad \Box \qquad \Box$

**Distributivity:** given $f(s,t), g(s,t), h(s,t)$ it results that

$(f * (g + h))(s,t) = (f * g)(s,t) + (f * h)(s,t)$. Left distributivity

$((f + g) * h)(s,t) = (f * h)(s,t) + (g * h)(s,t)$. Right distributivity

**Proof** Indeed, it results:

$$(f * (g + h))(s,t) = \int_s^t (g + h)(s,\tau) f(\tau,t) d\tau = \int_s^t (g(s,\tau) f(\tau,t) + h(s,\tau) f(\tau,t)) d\tau$$

$$= (f * g)(s,t) + (f * h)(s,t) \qquad \Box$$

**Bi-linearity:** The additivity is ensured by the distributivity. It remains to be proven that $\forall a \in \mathbb{R}, \forall f(s,t), g(s,t)$ it results:

$$(a(f * g))(s,t) = ((af) * g)(s,t) = (f * (ag))(s,t).$$

**Proof**

$$(a(f * g)(s,t)) = a\int_s^t g(s,\tau)f(\tau,t)d\tau = \int_s^t g(s,\tau)af(\tau,t)d\tau = ((af) * g)(s,t)$$

$$a\int_s^t g(s,\tau)f(\tau,t)d\tau = \int_s^t ag(s,\tau)f(\tau,t)d\tau = (f * (a g))(s,t).$$ □

**Remark 2.2** The non-homogeneous convolution is non-commutative, i.e.:

$$\exists f(s,t), g(s,t): (f * g)(s,t) \neq (g * f)(s,t).$$ □

**Example 2.1** $f(s,t) = e^{3s+4t}$, $g(s,t) = e^{-4s+2t}$ and their convolution is not commutative.

**Remark 2.3.** Time homogeneity is a particular case of non-homogeneity. Indeed, a two time variables case $f(s,t), 0 \leq s \leq t$ is homogeneous if

$$f(s,t) = f(s+h, t+h); \quad \forall 0 \leq s \leq t \text{ and } h: -s < h.$$

Given this hypothesis, it is possible to define another one-time variable function

$$\overline{f}(\tau) = f(s,t), \quad \forall s,t \text{ and } \tau = t - s. \tag{1}$$

It results that $\overline{f}$ is a function of the duration and not of the starting and arriving times. □

**Remark 2.4.** Given $f$ and $g$, the function of two time variables that are homogeneous in time. Let

$$\overline{f}(r) = f(t-s) \text{ and } \overline{g}(r) = g(t-s),$$

then from the homogeneity hypothesis and posed $\tau - s = x$ it results:

$$(f * g)(s,t) = \int_s^t g(s,\tau)f(\tau,t)d\tau = \int_s^t \overline{g}(\tau-s)\overline{f}(t-\tau)d\tau$$

$$= \int_0^u \overline{g}(x)\overline{f}(u-x)dx = (\overline{f} * \overline{g})(u), \text{ where } u = t - s$$
□

**3. Main definitions**

In this part, the main definitions of renewal theory, as given in Janssen & Manca (2006), are generalized in a non-homogeneous environment. Let $(X_n, n \geq 1)$ be a sequence of positive, independent and not necessarily identically distributed r.v. defined on the probability space $(\Omega, F, P)$.

**Definition 3.1** The random sequence $(T_n, n \geq 0)$, where:

$$T_0 = 0, \text{ and } T_n = X_1 + \cdots + X_n, \quad n \geq 1,$$

is called a *generalized renewal sequence* or sometimes *general renewal process*. □

The r.v. $T_n, n \geq 0$ are called *renewal times* and the r.v. $X_n, n \geq 1$ are called *interarrival times*. To avoid triviality it is supposed that $P[X_n \leq 0 = 0]$.

**Definition 3.4** The random sequence $(T_n, n \geq 0)$, is an *homogeneous classical renewal process* if the $X_n, n \geq 1$ r.v. are identically distributed and independent with $F$ as d.f.. □

**Defintion 3.5** The random sequence $(T_n, n \geq 0)$, is *a non-homogeneous classical renewal process* if variables $X_n, n \geq 1$ are conditionally independent (see Dawid (1979)). Given

$$(T_n, n \geq 1) \text{ where}$$
$$T_n = \sum_{i=1}^{n} X_i, n \geq 1$$

with here the function $F$ is the conditionally d.f. of $X_n$ given $T_n$ defined as follows

$$F(s,t) = P(X_n \leq t - s | T_{n-1} = s), 0 < s \leq t$$

And for $s = 0$, we have

$$F(0,t) = P(X_1 \leq t), 0 \leq t.$$

It follows that the process $(T_n, n \geq 0)$ is a non-homogenous Markov process with values in $\mathbb{R}^+$ □

**Remark 3.2**. In delayed homogeneous renewal processes $X_1$ has a different distribution respect to $X_n, n \geq 2$ □

**Definition 3.2.** With each renewal sequence, the following stochastic process can be associated with values in $\mathbb{N}$:

$$\left( N(t), \quad t \in \mathbb{R}^+ \right), \tag{2}$$

where as in the homogeneous and non-homogeneous cases, it respectively results:

$$\begin{aligned} N(t) > n-1 &\Leftrightarrow T_n \leq t, \quad n \in \mathbb{N} \\ N(t) - N(s) > n-1 &\Leftrightarrow T_{n'} = s, T_{n''} \leq t, n''-n' > n-1. \end{aligned} \tag{3}$$

These processes are respectively called the *homogenous and non-homogeneous associated counting process* or the *homogenous and non-homogeneous renewal counting process*. □

$N(t)$ represents the total number of "renewals" on $(0, t]$.

**Remark 3.3.** The probability of having at least $n$ renewals within a time $t$ in the homogeneous case and from time $s$ to time $t$ in the non-homogeneous environment is respectively given by:

$$\begin{aligned} \mathrm{P}[N(t) > n-1] &= F^{(n)}(t) \\ \mathrm{P}[N(t) - N(s) > n-1] &= F^{(n)}(s,t). \end{aligned}$$

The probability of having just $n$ renewals is obtained in the following way:

$$\mathrm{P}[N(t) = n] = \mathrm{P}[N(t) > n-1] - \mathrm{P}[N(t) > n] = F^{(n)}(t) - F^{(n+1)}(t),$$

$$\mathrm{P}[N(t) - N(s) = n] = \mathrm{P}[N(t) - N(s) > n-1] - \mathrm{P}[N(t) - N(s) > n] = F^{(n)}(s,t) - F^{(n+1)}(s,t). \quad \square$$

**Definition 3.3.** The homogeneous and non-homogeneous *renewal functions* are defined respectively as:

$$H(t) = E[N(t)], \ 0 \le t, t \in \mathbb{R}_0^+,$$
$$H(s,t) = E[N(t) - N(s)] = E[N(t)] - E[N(s)], \ 0 \le s \le t \ s,t \in \mathbb{R}_0^+,$$
(4)

provided that the expectation is finite. They give the mean number of renewals that verified respectively within a time $t$ and from time $s$ to time $t$. □

**Definition 3.4**. Given two non-homogeneous distribution functions (n.h.d.f.) $F(s,t)$ and $G(s,t)$, their convolution operation is defined in the following way:

$$G * F(s,t) = \int_s^t G(\tau,t) dF(s,\tau).$$ □

**Proposition 3.1.** The homogeneous and non-homogeneous *continuous time evolution equation* of the renewal equations, supposing the absolute continuity of the d.f. $F$, i.e. $dF(t) = f(t)dt$, can be written respectively in the following way:

$$H(t) = F(t) + \int_0^t f(\tau) H(t-\tau) d\tau,$$
$$H(s,t) = F(s,t) + \int_s^t f(s,\tau) H(\tau,t) d\tau.$$
(5)

**Proof.** For the homogeneous case see, for example, Janssen & Manca (2006).

In the non-homogeneous case it results:

$$E[N(t) - N(s)] = H(s,t) = \sum_{n=1}^{\infty} n \left[ F^{(n)}(s,t) - F^{(n+1)}(s,t) \right]$$

$$H(s,t) = \sum_{n=1}^{\infty} F^{(n)}(s,t) = F(s,t) + \sum_{n=2}^{\infty} F^{(n)}(s,t)$$
$$= F(s,t) + \left( \sum_{n=2}^{\infty} F^{(n-1)} \right) * F(s,t) = F(s,t) + \left( \sum_{n=1}^{\infty} F^{(n)} \right) * F(s,t)$$
$$H(s,t) = F(s,t) + H * F(s,t)$$

$$H(s,t) = F(s,t) + \int_s^t f(s,\tau) H(\tau,t) d\tau$$ □

**Definition 3.5.** The homogeneous and non-homogeneous *discrete time evolution equations* of the renewal processes are respectively the following:

$$H(t) = F(t) + \sum_{x=1}^{t} v(\tau)H(t-\tau)$$
$$H(s,t) = F(s,t) + \sum_{x=s+1}^{t} v(s,\tau)H(\tau,t) \quad (6)$$

where $F(t)$ and $F(s,t)$ are respectively the waiting time for the homogenous and non-homogeneous d.f.. They give the probability of having a renewal within the time $t$, given that the system was followed from time 0 in the homogeneous case and from time $s$ in the non-homogeneous environment. Furthermore, it results that:

$$v(t) = \begin{cases} 0 & \text{if } t = 0 \\ F(t) - F(t-1) & \text{if } t > 0 \end{cases}$$

and

$$v(s,t) = \begin{cases} 0 & \text{if } s \geq t \\ F(s,t) - F(s,t-1) & \text{if } s < t. \end{cases} \qquad \square$$

**Remark 3.4.** In the homogeneous case, the definitions of continuous and discrete time evolution equations correspond to those given respectively in Feller (1971) page 185 and Feller (1968) page 332. $\square$

## 4. Solution of non-homogeneous discrete time evolution equation

The discrete time non-homogeneous renewal equations (6) can be numerically solved very easily (for homogeneous case see Janssen & Manca (2006)). Compactly expressed, the non-homogeneous case of relation (6) can be written as:

$$H(s,t) - \sum_{x=s+1}^{t} v(s,\tau)H(\tau,t) = F(s,t) \Leftrightarrow \mathbf{U} \cdot \mathbf{H} = \mathbf{F} \quad (7)$$

where:

$$\mathbf{U} = \begin{bmatrix} 1 & -v_{0,1} & \cdots & -v_{0,k} & \cdots & -v_{0,n} & \cdots \\ 0 & 1 & \cdots & -v_{1,k} & \cdots & -v_{1,n} & \cdots \\ \vdots & \vdots & \ddots & \vdots & \ddots & \vdots & \ddots \\ 0 & 0 & \cdots & 1 & \cdots & -v_{k,n} & \cdots \\ \vdots & \vdots & \ddots & \vdots & \ddots & \vdots & \ddots \\ 0 & 0 & \cdots & 0 & \cdots & 1 & \cdots \\ \vdots & \vdots & \ddots & \vdots & \ddots & \vdots & \ddots \end{bmatrix},$$

is the coefficient matrix;

$$\mathbf{H} = \begin{bmatrix} H_{0,0} & H_{0,1} & \cdots & H_{0,k} & \cdots & H_{0,n} & \cdots \\ 0 & H_{1,1} & \cdots & H_{1,k} & \cdots & H_{1,n} & \cdots \\ \vdots & \vdots & \ddots & \vdots & \ddots & \vdots & \ddots \\ 0 & 0 & \cdots & H_{k,k} & \cdots & H_{k,n} & \cdots \\ \vdots & \vdots & \ddots & \vdots & \ddots & \vdots & \ddots \\ 0 & 0 & \cdots & 0 & \cdots & H_{n,n} & \cdots \\ \vdots & \vdots & \ddots & \vdots & \ddots & \vdots & \ddots \end{bmatrix},$$

is the unknown matrix and

$$\mathbf{F} = \begin{bmatrix} F_{0,0} & F_{0,1} & \cdots & F_{0,k} & \cdots & F_{0,n} & \cdots \\ 0 & F_{1,1} & \cdots & F_{1,k} & \cdots & F_{1,n} & \cdots \\ \vdots & \vdots & \ddots & \vdots & \ddots & \vdots & \ddots \\ 0 & 0 & \cdots & F_{k,k} & \cdots & F_{k,n} & \cdots \\ \vdots & \vdots & \ddots & \vdots & \ddots & \vdots & \ddots \\ 0 & 0 & \cdots & 0 & \cdots & F_{n,n} & \cdots \\ \vdots & \vdots & \ddots & \vdots & \ddots & \vdots & \ddots \end{bmatrix},$$

is the matrix in which each row is a discrete time d.f..

**Remark 4.1.** System (7) allows for the solution. Indeed, the determinant of the coefficient matrix of the system (7) is equal to 1 (see Riesz (1913)). □

**Remark 4.2.** $F(k,k)=0$ because it is impossible to have a renewal in a time 0. This implies $v(k,k)=0$ and $H(k,k)=0$. □

### 4.1. Some particular formulae

In this part, some formulas of the numerical solution of the non-homogeneous part of renewal equations (5) will be given. The relations will be related to particular generalized Newton-Cotes formulas (see Hildebrand (1987)).

Simpson quadrature method gives the following equations:

$$\hat{H}(uh,kh) = F(uh,kh) + \frac{h}{3}\hat{H}(uh,kh)f(uh,uh) + \frac{4h}{3}\sum_{\tau=\left\lfloor\frac{u}{2}\right\rfloor+1}^{\left\lfloor\frac{k}{2}\right\rfloor} \hat{H}((2\tau-1)h,kh)f(uh,(2\tau-1)h)$$

$$+\frac{2h}{3}\sum_{\tau=\left\lfloor\frac{u}{2}\right\rfloor+1}^{\left\lfloor\frac{k}{2}\right\rfloor-1} \hat{H}(2\tau h,kh)f(uh,2\tau h) + \frac{h}{3}\hat{H}(kh,kh)f(uh,kh),$$

where $\hat{H}$ represents the approximate value of $H$ and $h$ the discretization interval.

Using Bezout's quadrature method the following relation is obtained:

$$\hat{H}(uh,kh) = F(uh,kh) + \frac{h}{2}\hat{H}(uh,kh)f(uh,uh) + h\sum_{\tau=u+1}^{k-1} \hat{H}(\tau h,kh)f(uh,\tau h)$$

$$+\frac{h}{2}\hat{H}(kh,kh)f(uh,kh).$$

Finally, if the simplest quadrature method (rectangle formula) is applied, then it is possible to obtain two different formulas; one giving the value of the integrating function at the end of the interval and the other at the beginning of the interval. In this way, the following formulas are obtained:

$$\hat{H}(uh,kh) = F(uh,kh) + h\sum_{\tau=u+1}^{k} \hat{H}(\tau h,kh)f(uh,\tau h).$$

$$\hat{H}(uh,kh) = F(uh,kh) + h\sum_{\tau=u}^{k-1} \hat{H}(\tau h,kh)f(uh,\tau h).$$

Substituting the differential $(hf(uh,kh))$ by means of the difference, it respectively results:

$$\hat{H}(uh,kh) \cong F(uh,kh) + \sum_{\tau=u+1}^{k} \hat{H}(\tau h, kh)(F(uh,\tau h) - F(uh,(\tau-1)h)),$$

$$\hat{H}(uh,kh) \cong F(uh,kh) + \sum_{\tau=u}^{k-1} \hat{H}(\tau h, kh)(F(uh,(\tau+1)h) - F(uh,(\tau)h)).$$

(8)

**Remark 4.3** In Baker (1977) (see page 925), there are two lemmas and a theorem ensuring that, under our conditions, the approximation of (8) to the solution tends to 0 as the discretization interval $h$ tends to 0. But, given our particular system, in the following we are able to give a more interesting result. □

### 4.2. Relations between discrete time and continuous time renewal equations

Posing $h = 1$ in (8)

$$\hat{H}(u,k) \cong F(u,k) + \sum_{\tau=u+1}^{k} \hat{H}(\tau,k)(F(u,\tau) - F(u,\tau-1)).$$

Furthermore, giving the following positions:

$$v(u,\tau) = \begin{cases} F(u,\tau) = 0 & \tau = u, \\ F(u,\tau) - F(u,\tau-1) & \tau > u, \end{cases}$$

the following is obtained:

$$H(u,k) = F(u,k) + \sum_{\tau=1}^{k} H(\tau,k) v(u,\tau),$$

(9)

that is the discrete time non-homogeneous renewal equation.

Now let $H$ be a continuous time renewal function and $\{T_n\}$ the related renewal process.

If the following is set:

$$T_n^h = \left\lfloor \frac{T_n}{h} \right\rfloor h$$

and

$$N^h(t) = n \text{ if } T_n^h \le t < T_{n+1}^h, \tag{10}$$

then the related discrete time renewal function is given by:

$$H^h(uh, kh) = F^h(uh, kh) + \sum_{\tau=u+1}^{k} v^h(uh, \tau h) H^h((k-\tau)h, kh).$$

The renewal process $T_n^h$ is defined in the same probability space $(\Omega, F, P)$ of $T_n$.

Given $\omega \in \Omega$ the following result holds:

**Proposition 4.1** *The $T_n^h$ process converges to $T_n$ for $h \to 0$.* □

**Proof.** From the definitions given in section 3 and (10) it results:

$$P\left[\underset{h \to 0}{N^h(t) = n}\right] \to P[N(t) = n] \quad \forall n \in \mathbb{N}, t \in \mathbb{R}^+. \tag{11}$$

The formula (11) implies:

$$T_n^h \xrightarrow[h \to 0]{a.s.} T_n.$$
□

**Remark 4.4.** The **Proposition 4.1** can be considered a very special case of the **Theorem 10** given in Janssen Manca (2001). But in this paper, the proof is trivial. □

**Remark 4.5.** The results obtained show that it is possible to obtain the non-homogeneous discrete time renewal equation by means of the simplest discretization of the related continuous time and that, starting from the discrete time, it is also possible to obtain the continuous time. □

**Remark 4.6.** Freiberger & Grenander (1971) demonstrate a similar approach to the renewal equation numerical solution for the homogeneous case, but there is no justification of the method. Many other papers (see the introduction) deal with the same problem in the homogeneous case, but, as far as the authors know, the relationship between the discrete time and continuous time renewal process has only been justified in the book by Janssen & Manca (2006), as they are in this paper in the non-homogeneous case. It is, however, the first time that the numerical treatment of the non-homogenous renewal processes is presented. ▫

## 5. Car insurance mean number claims calculation in a non-homogeneous age environment

**5.1. The construction of homogeneous and non-homogeneous d.f.**

In this section, we describe the construction of d.f. from the raw data. We will show firstly, the homogeneous case and secondly, the non-homogeneous.

Each insured individual corresponds to a record. For each record, the necessary data in the homogeneous case are:

- how many accidents the subject caused in the observed time horizon and
- the date of each accident.

The vector of the number of claims that will be paid from the insurance company will be constructed from these data. In $i^{th}$ elements of the vector, there will be the number of accidents that occurred just in $i$ years from the first car insurance contract or from the previous accident. At the end of this process, the vector will have the following shape:

| $n(1)$ | $n(2)$ | … | $n(i)$ | … | $n(T)$ |
|---|---|---|---|---|---|

**Table 1:** Number of claims within the time $i$

$T$ represents the horizon time and $T = \max(i)$ where $i$ is all the possible times obtained during the observation of the data. Subsequently, another vector will be created whose elements are defined in the following way:

$$v(i) = \sum_{k=1}^{i} n(i)$$

The next step will create the d.f. i.e.:

$$F(i) = \frac{v(i)}{v(T)}, \ \forall i = 1,\ldots,T.$$

In the non-homogeneous case, a vector will be created for each age $s$, i.e.:

|  | Arriving age | | | | | | |
|---|---|---|---|---|---|---|---|
| Starting age | $n(0,1)$ | $n(0,2)$ | $n(0,3)$ | ... | $n(0,t)$ | ... | $n(0,T)$ |
| | 0 | $n(1,2)$ | $n(1,3)$ | ... | $n(1,t)$ | ... | $n(1,T)$ |
| | 0 | 0 | $n(2,3)$ | ... | $n(2,t)$ | ... | $n(2,T)$ |
| | ⋮ | ⋮ | ⋮ | ⋱ | ⋮ | ⋱ | ⋮ |
| | 0 | 0 | 0 | ... | $n(t-1,t)$ | ... | $n(t-1,T)$ |
| | ⋮ | ⋮ | ⋮ | ⋱ | ⋮ | ⋱ | ⋮ |
| | 0 | 0 | 0 | ... | 0 | ... | $n(T-1,T)$ |

**Table 2:** Number of claims reported starting from age $s$ just at age $t>s$

**Remark 5.1.** The age 0 corresponds to the starting age 18. The d.f. are constructed for each starting age. At time 0 it is impossible to have a claim, for this reason in the main diagonal of Table 2 there are the *ages s*-1, *s*. The non-homogeneous d.f. $F(s,t)$ should be calculated as $\forall s = 0,\ldots,T-1$ in the same way as the $F(t)$ in the homogeneous case.

## 5.2. The data description

In this section, the non-homogeneous renewal equation is applied to a real actuarial data application. It is demonstrated that the renewal process can be applied in a very general case utilising data from statistical observations. The relation that will be used in this case is (9).

The renewal process will be applied in the case of the calculation of the mean number of claims; common problem for motorcar insurance companies. In this insurance, as it is well known, the age of the insured person assumes great relevance. It is possible to take into account the different behaviours of insured people as a function of age by means of a non-homogeneous renewal process where the non-homogeneity is related to age.

In a motorcar insurance contract, each time the insured has an accident the insurance company will pay for the damage. In our model, the renewed contract will consider the driver's age. Our non-homogeneous environment takes into account this fact. Indeed the d.f. were different as a function of the starting age, but the independence hypothesis holds.

We have raw data regarding accidents that an insurance company collected during a period of about 50 years up to the year 2000. It is possible to construct the non-homogeneous discrete time d.f. of the renewal time of the motorcar claims from this data. The data covered a total of 156,428 insured people of which 22,395 had had at least one accident. Data concerning the insurance premiums were not available. The dates of individual contracts were available and 60,278 of those appeared to be correct. The ages of all the insured people were known. We computed the mean age upon entering the contract of these 60,278 records and we found that the mean age was 23.77 years.

We supposed that all the other 96,150 insured people entered the contract at age 24.

We had 32,201 claims in total. At most, we had three claims per insured person. For example, the renewal worked in this way; it was supposed that an insured person entered her/his contract at 23 years. She/he made the first claim at 41. We took into account that, for this insured person of starting age 23, there was a claim after 18 years. The renewed age was now 41. It was supposed that she/he made another claim at age 50. In this way, we also suppose that a person that was 41 years old made a new claim after 9 years.

Before applying our model, we did a preliminary study on the data working in a homogeneous environment. To get correct results, we worked on the second and third accidents because, in this way, we were sure about the time of the renewal. The dates of claims were considered to be accurate. The results are reported in the **Table 5.3**.

In this table, in the second and third column the number of accidents are reported after 1 year, then after 2 years and so on, are given. The heading 1-2 means that the number of claims is related to the second accident and that the waiting time is given by the time that passed between the first and the second claim. This waiting time is reported in the first column. The heading of the third column has the same meaning and its contents report the number of claims that occurred at the waiting time given in the first column. For example, 1,576 in the third row of this table represents the number of second accidents that occurred after two years and before the third year since the first accident. The number 226 represents the number of third accidents that occurred after two years and before the third year since the second accident. The fourth and fifth columns give the probability function of waiting times that are related to the columns 2 and 3, respectively.

| # years | 1-2 | 2-3 | prob 1-2 | prob 2-3 |
|---|---|---|---|---|
| 1 | 153 | 88 | 0.018595 | 0.055767 |
| 2 | 695 | 126 | 0.084468 | 0.079848 |
| 3 | 1576 | 226 | 0.191541 | 0.143219 |
| 4 | 1549 | 224 | 0.18826 | 0.141952 |
| 5 | 1344 | 172 | 0.163345 | 0.108999 |
| 6 | 928 | 176 | 0.112786 | 0.111534 |
| 7 | 619 | 138 | 0.075231 | 0.087452 |
| 8 | 386 | 107 | 0.046913 | 0.067807 |
| 9 | 278 | 86 | 0.033787 | 0.054499 |
| 10 | 189 | 69 | 0.02297 | 0.043726 |
| 11 | 139 | 61 | 0.016894 | 0.038657 |
| 12 | 101 | 33 | 0.012275 | 0.020913 |
| 13 | 77 | 21 | 0.009358 | 0.013308 |
| 14 | 40 | 15 | 0.004861 | 0.009506 |
| 15 | 36 | 16 | 0.004375 | 0.010139 |
| 16 | 15 | 7 | 0.001823 | 0.004436 |
| 17 | 11 | 3 | 0.001337 | 0.001901 |
| 18 | 11 | 4 | 0.001337 | 0.002535 |
| 19 | 10 | 3 | 0.001215 | 0.001901 |
| 20 | 7 | 2 | 0.000851 | 0.001267 |
| 21 | 14 | 1 | 0.001702 | 0.000634 |
| 22 | 6 | 0 | 0.000729 | 0 |
| 23 | 8 | 0 | 0.000972 | 0 |
| 24 | 7 | 0 | 0.000851 | 0 |
| 25 | 7 | 0 | 0.000851 | 0 |
| 26 | 2 | 0 | 0.000243 | 0 |
| 27 | 0 | 0 | 0 | 0 |
| 28 | 4 | 0 | 0.000486 | 0 |
| 29 | 3 | 0 | 0.000365 | 0 |
| 30 | 5 | 0 | 0.000608 | 0 |
| 31 | 3 | 0 | 0.000365 | 0 |
| 32 | 2 | 0 | 0.000243 | 0 |
| 33 | 1 | 0 | 0.000122 | 0 |
| 34 | 0 | 0 | 0 | 0 |
| 35 | 1 | 0 | 0.000122 | 0 |
| 36 | 0 | 0 | 0 | 0 |
| 37 | 1 | 0 | 0.000122 | 0 |
| Total | 8228 | 1578 | 1 | 1 |

**Table 3: Homogeneous study of II and III claim**

**Figure 1** reports the histograms of the last two columns of **Table 3**.

From the first study of data we knew:

1 - 14.3% of our population had accidents,

2 - less than 40% of the information is totally reliable.

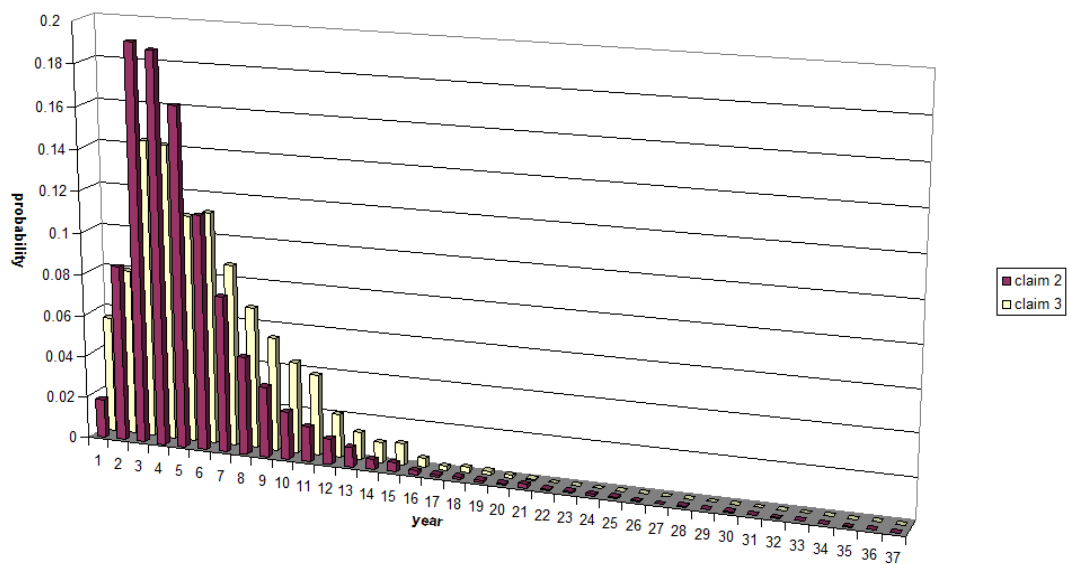

**Figure 1: Probability distribution obtained by data of the II and III claim**

It is clear from the shape of the histograms that the discrete time probability distribution is not geometric and this implies that our renewal process does not abide by Poisson's function.

| age | tot numb | no claim | prob no claim | prob claim |
|---|---|---|---|---|
| 18 | 279 | 237 | 0.849462 | 0.150538 |
| 19 | 8863 | 6745 | 0.761029 | 0.238971 |
| 20 | 12960 | 9893 | 0.763349 | 0.236651 |
| 21 | 11745 | 9346 | 0.795743 | 0.204257 |
| 22 | 4388 | 3341 | 0.761395 | 0.238605 |
| 23 | 3137 | 2335 | 0.744342 | 0.255658 |
| 24 | 2100 | 1580 | 0.752381 | 0.247619 |
| 25 | 1651 | 1248 | 0.755906 | 0.244094 |
| 26 | 1397 | 1058 | 0.757337 | 0.242663 |
| 27 | 1092 | 796 | 0.728938 | 0.271062 |
| 28 | 997 | 748 | 0.750251 | 0.249749 |
| 29 | 1037 | 748 | 0.721311 | 0.278689 |
| 30 | 944 | 678 | 0.718220 | 0.28178 |
| 31 | 1120 | 878 | 0.783929 | 0.216071 |
| 32 | 722 | 534 | 0.739612 | 0.260388 |
| 33 | 659 | 502 | 0.761760 | 0.23824 |
| 34 | 555 | 425 | 0.765766 | 0.234234 |
| 35 | 486 | 372 | 0.765432 | 0.234568 |
| 36 | 436 | 329 | 0.754587 | 0.245413 |
| 37 | 408 | 303 | 0.742647 | 0.257353 |
| 38 | 371 | 274 | 0.738544 | 0.261456 |
| 39 | 380 | 284 | 0.747368 | 0.252632 |
| 40 | 412 | 291 | 0.70631 | 0.293689 |
| 41 | 461 | 339 | 0.735358 | 0.264642 |
| 42 | 313 | 240 | 0.766773 | 0.233227 |
| 43 | 289 | 207 | 0.716263 | 0.283737 |
| 44 | 283 | 212 | 0.749117 | 0.250883 |
| 45 | 248 | 181 | 0.72984 | 0.270161 |
| 46 | 253 | 198 | 0.782609 | 0.217391 |
| 47 | 214 | 156 | 0.728972 | 0.271028 |
| 48 | 200 | 144 | 0.72 | 0.28 |
| 49 | 224 | 159 | 0.709821 | 0.290179 |
| 50 | 181 | 135 | 0.745856 | 0.254144 |
| 51 | 168 | 141 | 0.839286 | 0.160714 |
| 52 | 141 | 107 | 0.758865 | 0.241135 |
| 53 | 127 | 107 | 0.842520 | 0.15748 |
| 54 | 121 | 106 | 0.876033 | 0.123967 |
| 55 | 100 | 81 | 0.81 | 0.19 |
| 56 | 99 | 84 | 0.848485 | 0.151515 |
| 57 | 87 | 68 | 0.781609 | 0.218391 |
| 58 | 66 | 57 | 0.863636 | 0.136364 |
| 59 | 70 | 63 | 0.9 | 0.1 |
| ≥60 | 600 | 535 | 0.891666667 | 0.108333 |
| total | 60384 | 46265 | 0.766179783 | 0.23382 |

**Table 4: Probability of not having claims**

Given that a non-homogeneous model had to be constructed, we decided to take into account all the data related to an accident. To do this, an entrance age of 24 years was given to all the people that did not have the correct age upon entering the contract. We did not consider the claims with a date earlier than the starting age of contract or of the previous renewal. In the end, we considered 23,395 claims and discarded 8,806 accidents.

We also constructed a statistic on the number of people that did not have claims as a function of their age upon entering the contract. The statistic was carried out on the 60,278 reliable records. We report these data in **Table 4**. In these data, which are more reliable than the data of the complete file, a higher probability of having accidents was obtained. The global probability of having an accident becomes 23.3% instead of 9%. Furthermore, we constructed the matrix of occurrences taking into account all the data, i.e. for each starting age from 18 to 60, (after 60 there were fewer new contracts so we compiled them all together), we constructed the number of claims that were made according to the age of the person who caused the accident, see **Figure 2**. Normalizing for each age, we obtained the probability function for each starting age. The shape of the probability function is the same as the shape of occurrences. In this case, too, it is evident that, as in the homogeneous case, the probability functions do not have the same shape as a geometric probability distribution and the Poisson hypotheses should be rejected.

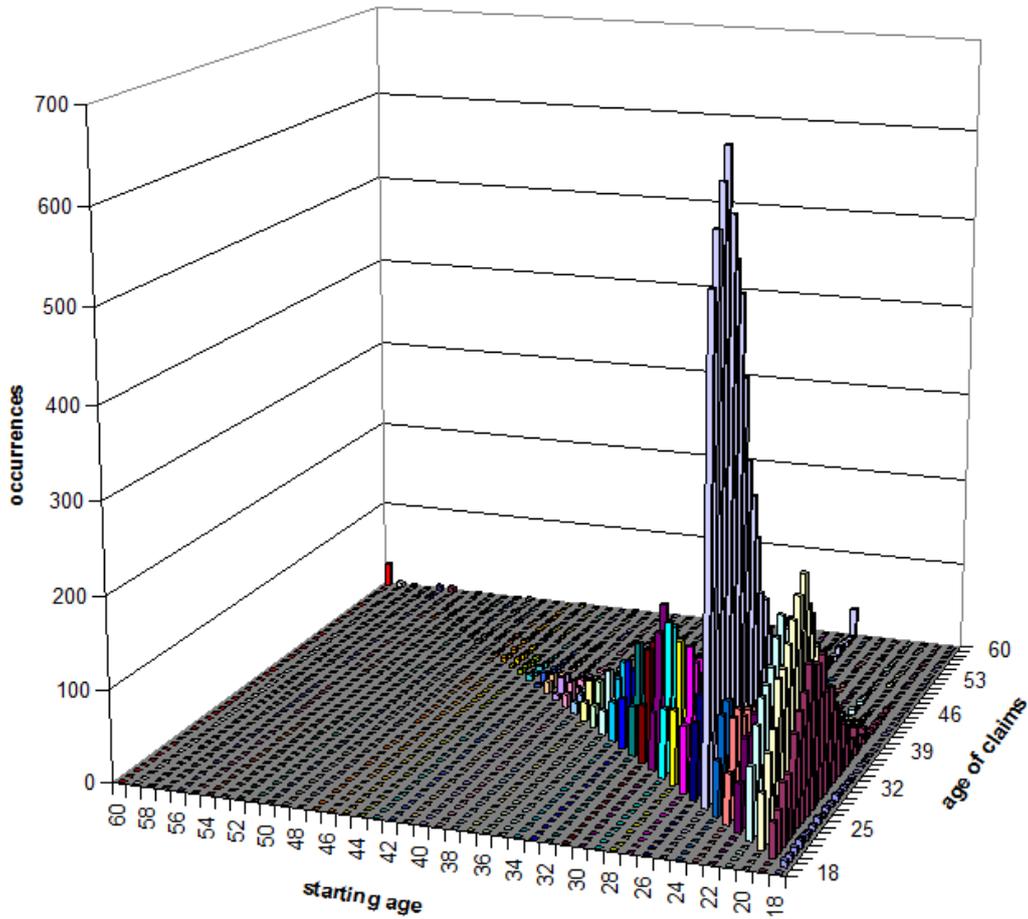

**Figure 2: Occurrence distribution of claims as a function of the contract age**

*5.3 The result description*

In light of these results, we decided to work with the physical measure. We constructed the waiting time probability functions for each starting age from the occurrences. From these probability functions we constructed the cumulative d.f. and we could then apply the relation (9).

To show just how simple it is to solve the non-homogeneous renewal equation and, in this way, to obtain the mean number of renewals for each considered age, we report the program written in Mathematica 9 language.

```
frip=Table[0.0,{i,1,43},{j,1,43}];
phi=Table[0.0,{i,1,43},{j,1,43}];
For[i=1,i≤nanni,i++,
  frip[[i]]=ReadList[puntf,Number,nanni];
];
For[h=nanni,h≥2,h--,
  For[k=h-1,k≥1,k--,
   phi[[k,h]]=frip[[k,h]];
   For[i=k+1, i≤h-1, i++,
    phi[[k,h]]+=phi[[i,h]]*(frip[[k,i]]-frip[[k,i-1]]);
   ];
  ];
];
```

In **Figure 3** the non-homogeneous d.f. obtained from the data are given.

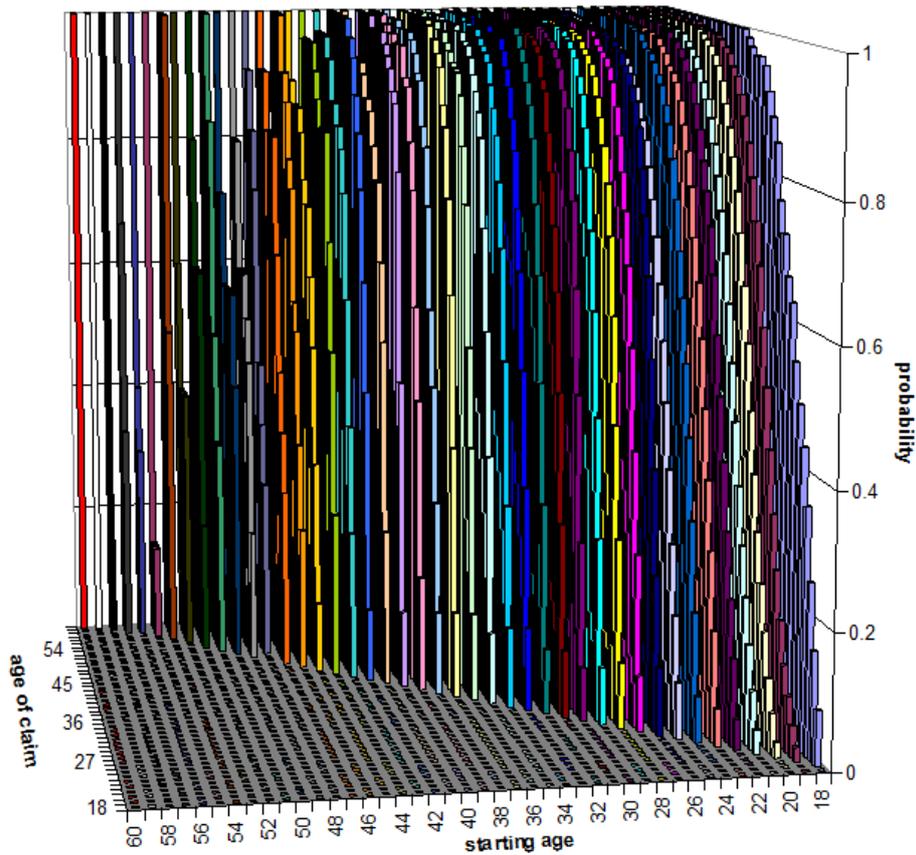

**Figure 3: Waiting time d.f. as a function of age**

Finally, in **Table 5**, the mean number of claims for starting ages 20 to 30 and for each age of claims is given. The comparison should be done along the diagonal because, in this way, the same number of years is considered.

**6. Conclusions**

In the authors' opinion, the main results of this paper are:

1 - the defining of non-homogeneous in time convolutions and their properties,

2 - the systematization of the theory of non-homogeneous renewal processes,

3 - the numerical treatment of non-homogeneous renewal processes,

4 - the application of the method to the motorcar insurance environment, and

5 - the possibility, within the renewal equation, of using the d.f. constructed directly from the observed data explaining how simple it is to start from the real data and to apply the model.

|  |  | Age of contract |  |  |  |  |  |  |  |  |  |  |
|---|---|---|---|---|---|---|---|---|---|---|---|---|
|  |  | 20 | 21 | 22 | 23 | 24 | 25 | 26 | 27 | 28 | 29 | 30 |
| Mean number of claims done up to the reported age | 20 | 0 | 0 | 0 | 0 | 0 | 0 | 0 | 0 | 0 | 0 | 0 |
|  | 21 | 0.023 | 0 | 0 | 0 | 0 | 0 | 0 | 0 | 0 | 0 | 0 |
|  | 22 | 0.060 | 0.035 | 0 | 0 | 0 | 0 | 0 | 0 | 0 | 0 | 0 |
|  | 23 | 0.108 | 0.088 | 0.047 | 0 | 0 | 0 | 0 | 0 | 0 | 0 | 0 |
|  | 24 | 0.177 | 0.155 | 0.131 | 0.056 | 0 | 0 | 0 | 0 | 0 | 0 | 0 |
|  | 25 | 0.255 | 0.233 | 0.222 | 0.167 | 0.069 | 0 | 0 | 0 | 0 | 0 | 0 |
|  | 26 | 0.352 | 0.326 | 0.340 | 0.296 | 0.192 | 0.084 | 0 | 0 | 0 | 0 | 0 |
|  | 27 | 0.468 | 0.442 | 0.473 | 0.425 | 0.340 | 0.186 | 0.101 | 0 | 0 | 0 | 0 |
|  | 28 | 0.604 | 0.578 | 0.613 | 0.570 | 0.493 | 0.305 | 0.240 | 0.084 | 0 | 0 | 0 |
|  | 29 | 0.771 | 0.730 | 0.777 | 0.741 | 0.663 | 0.452 | 0.423 | 0.266 | 0.098 | 0 | 0 |
|  | 30 | 0.945 | 0.904 | 0.959 | 0.915 | 0.845 | 0.609 | 0.636 | 0.446 | 0.285 | 0.087 | 0 |
|  | 31 | 1.133 | 1.083 | 1.144 | 1.094 | 1.023 | 0.779 | 0.842 | 0.655 | 0.482 | 0.280 | 0.081 |
|  | 32 | 1.333 | 1.278 | 1.336 | 1.307 | 1.233 | 0.969 | 1.053 | 0.870 | 0.709 | 0.491 | 0.272 |
|  | 33 | 1.542 | 1.482 | 1.543 | 1.519 | 1.450 | 1.168 | 1.274 | 1.083 | 0.938 | 0.719 | 0.520 |
|  | 34 | 1.768 | 1.704 | 1.770 | 1.742 | 1.682 | 1.378 | 1.514 | 1.305 | 1.167 | 0.958 | 0.761 |
|  | 35 | 1.995 | 1.930 | 2.000 | 1.971 | 1.905 | 1.596 | 1.746 | 1.539 | 1.397 | 1.200 | 1.013 |
|  | 36 | 2.218 | 2.142 | 2.220 | 2.189 | 2.125 | 1.804 | 1.961 | 1.760 | 1.613 | 1.422 | 1.236 |
|  | 37 | 2.474 | 2.394 | 2.469 | 2.437 | 2.374 | 2.047 | 2.216 | 2.024 | 1.867 | 1.688 | 1.505 |
|  | 38 | 2.724 | 2.642 | 2.721 | 2.681 | 2.625 | 2.289 | 2.462 | 2.273 | 2.122 | 1.939 | 1.750 |
|  | 39 | 2.975 | 2.892 | 2.971 | 2.933 | 2.874 | 2.530 | 2.718 | 2.523 | 2.371 | 2.196 | 2.004 |
|  | 40 | 3.247 | 3.162 | 3.242 | 3.202 | 3.141 | 2.794 | 2.985 | 2.793 | 2.641 | 2.471 | 2.280 |
|  | 41 | 3.500 | 3.412 | 3.488 | 3.454 | 3.391 | 3.038 | 3.238 | 3.042 | 2.889 | 2.724 | 2.529 |
|  | 42 | 3.799 | 3.712 | 3.790 | 3.749 | 3.692 | 3.329 | 3.532 | 3.344 | 3.189 | 3.021 | 2.830 |
|  | 43 | 4.118 | 4.028 | 4.108 | 4.067 | 4.008 | 3.643 | 3.850 | 3.665 | 3.507 | 3.344 | 3.146 |
|  | 44 | 4.420 | 4.329 | 4.410 | 4.366 | 4.312 | 3.939 | 4.152 | 3.969 | 3.811 | 3.648 | 3.453 |
|  | 45 | 4.675 | 4.584 | 4.666 | 4.621 | 4.568 | 4.188 | 4.407 | 4.222 | 4.067 | 3.903 | 3.708 |
|  | 46 | 4.924 | 4.832 | 4.915 | 4.867 | 4.818 | 4.435 | 4.654 | 4.472 | 4.316 | 4.153 | 3.960 |
|  | 47 | 5.200 | 5.109 | 5.191 | 5.143 | 5.095 | 4.709 | 4.929 | 4.748 | 4.591 | 4.430 | 4.240 |
|  | 48 | 5.476 | 5.386 | 5.468 | 5.418 | 5.371 | 4.982 | 5.205 | 5.025 | 4.868 | 4.706 | 4.518 |
|  | 49 | 5.617 | 5.527 | 5.609 | 5.558 | 5.513 | 5.122 | 5.345 | 5.166 | 5.009 | 4.849 | 4.660 |
|  | 50 | 5.722 | 5.632 | 5.714 | 5.665 | 5.617 | 5.226 | 5.449 | 5.271 | 5.114 | 4.956 | 4.764 |
|  | 51 | 5.999 | 5.909 | 5.991 | 5.940 | 5.894 | 5.500 | 5.723 | 5.548 | 5.390 | 5.232 | 5.043 |
|  | 52 | 6.239 | 6.148 | 6.230 | 6.179 | 6.133 | 5.738 | 5.961 | 5.787 | 5.631 | 5.472 | 5.284 |
|  | 53 | 6.398 | 6.307 | 6.389 | 6.339 | 6.292 | 5.897 | 6.119 | 5.946 | 5.791 | 5.631 | 5.443 |
|  | 54 | 6.578 | 6.487 | 6.569 | 6.519 | 6.472 | 6.075 | 6.298 | 6.127 | 5.971 | 5.811 | 5.624 |
|  | 55 | 6.733 | 6.644 | 6.725 | 6.675 | 6.628 | 6.230 | 6.454 | 6.283 | 6.126 | 5.967 | 5.779 |
|  | 56 | 6.965 | 6.875 | 6.957 | 6.906 | 6.859 | 6.461 | 6.687 | 6.515 | 6.358 | 6.197 | 6.011 |
|  | 57 | 7.131 | 7.041 | 7.123 | 7.072 | 7.025 | 6.627 | 6.853 | 6.681 | 6.524 | 6.365 | 6.178 |
|  | 58 | 7.400 | 7.310 | 7.392 | 7.340 | 7.293 | 6.895 | 7.121 | 6.950 | 6.792 | 6.633 | 6.447 |
|  | 59 | 7.573 | 7.483 | 7.565 | 7.513 | 7.466 | 7.068 | 7.295 | 7.122 | 6.965 | 6.807 | 6.620 |
|  | 60 | 8.573 | 8.483 | 8.565 | 8.513 | 8.466 | 8.068 | 8.295 | 8.122 | 7.965 | 7.807 | 7.620 |

**Table 5: Mean number of claims – starting ages 20-30**

The proposed application in the field of automobile insurance could be interesting for insurance companies. Indeed, by being able to predict the mean number of claims a person will make on their motorcar insurance contract at a given age during his/her driving life, companies can judge as to whether an insured person is a good or a bad driver.

In the near future, we will attempt to also introduce the running time in our model. The construction of such a model will imply that the results will be a function of two time variables.